\documentclass[a4paper,twoside,10pt]{report}
\usepackage{amssymb, amsmath, amsfonts}
\usepackage{epsfig,amsfonts,amsbsy,amssymb,amsthm}
\usepackage{float,color}
\usepackage{fancyhdr}
\usepackage{makeidx}
\usepackage{graphics}
\usepackage{rupam-DE}
\usepackage{graphicx}
\usepackage{amsbsy}
\usepackage{amssymb}
\usepackage{amsmath}
\usepackage{graphics}
\usepackage{makeidx}
\usepackage{polynom}

\newcommand{\bea}{\begin{eqnarray}}
\newcommand{\eea}{\end{eqnarray}}
\newcommand{\ben}{\begin{eqnarray*}}
\newcommand{\een}{\end{eqnarray*}}

\theoremstyle{plain}
\newtheorem{thm}{Theorem}[section]
\newtheorem{cor}[thm]{Corollary}
\newtheorem{lem}[thm]{Lemma}
\newtheorem{prop}[thm]{Proposition}
\newtheorem{property}[thm]{Property}

\theoremstyle{definition}
\newtheorem{definition}{Definition}[section]
\newtheorem{exm}{Example}[section]

\newtheorem{ex}{Exercise}[section]
\newtheorem{note}{Note}[section]
\newtheorem{proposition}{Proposition}[section]
\theoremstyle{remark}
\newtheorem{rem}{Remark}[section]
\numberwithin{equation}{section}
\renewcommand{\theequation}{\thesection.\arabic{equation}}
\def\part{{\partial}}

\catcode`@=11 \@addtoreset{equation}{section}
\renewcommand\theequation{\thesection.\@arabic\c@equation}
\catcode`@=12

\makeatother
\makeindex

\begin{document}

\pagenumbering{roman}

\def \noin{\noindent}
\fancyfoot[CE,CO]{
{\thepage}}
\renewcommand{\footrulewidth}{0pt}
\renewcommand{\headrulewidth}{0.0pt}
\pagestyle{fancy}

\fancypagestyle{plain}{
\fancyhf{}
\fancyfoot[CO,CE]{
{\hspace{\textwidth}{}
{\thepage}}
}
\renewcommand{\headrulewidth}{0pt}
}

\makeatletter
\def\cleardoublepage{\clearpage\if@twoside \ifodd\c@page\else
	\hbox{}
	\vspace*{\fill}
	\thispagestyle{empty}
	\newpage
	\if@twocolumn\hbox{}\newpage\fi\fi\fi}
\makeatother


\baselineskip=22pt
\setcounter{page}{0}
{
\thispagestyle{empty}
\begin{center}
\vspace{1cm}

{\Huge\sf Gonit Sora Notes}\\
\vspace{1cm}
{\Huge\sf Probability Theory}\\

\vspace{3cm}

{\Large \sf By} \\

\vspace{2cm}
{\Large \sf Manjil P. Saikia} \\
{\Large \sf  Rupam Haloi}\\
{\large \sf \today}

\vspace{2cm}
{\Large \sf Department of Mathematical Sciences}\\
{\Large \sf Tezpur University}\\
{\Large \sf Napaam-784028, India}

\end{center}


\clearpage%
\begin{center}
\huge{\textbf{Preface}}
\end{center}

This is the first of the proposed sets of notes to be published in the website Gonit Sora (http://gonitsora.com). The notes will hopefully be able to help the students to learn their subject in an easy and comprehensible way. These notes are aimed at mimicking exactly what would be typically taught in a one-semester course at a college or university. The level of the notes would be roughly at the undergraduate level.

The present sets of notes are not yet complete and this is the second version that is being posted. These notes contain very few proofs and only state the important results in \textit{Probability Theory}. These notes are based on the course taught at \texttt{Tezpur University, Assam, India} by \textit{Dr. Santanu Dutta}. There may be some errors and typos in these notes which we hope the reader would bring to our notice.

\bigskip

\textbf{Manjil P. Saikia}

\smallskip

\textit{manjil@gonitsora.com}

\smallskip

\textbf{Rupam Haloi}

\baselineskip=20pt

\clearpage

\baselineskip=20pt


\baselineskip=20pt
\tableofcontents
\fancyhead[LE]{{Contents}}
\fancyhead[LO]{{\large\sffamily Contents}}
\renewcommand{\headrulewidth}{0.5pt}

\clearpage

\fancyhead[LO]{{\it Probability Theory: Unit \thechapter}}
\fancyhead[RO]{\it\thesection}
\fancyhead[RO]{\it\rightmark}

\fancyhead[RE]{{\it Probability Theory: Unit \thechapter}}
\fancyhead[LE]{\it\rightmark}
\renewcommand{\headrulewidth}{0.5pt}


\pagenumbering{arabic}
\pagestyle{fancy}

\setcounter{chapter}{0}
\setcounter{page}{1}

\baselineskip=20pt
\renewcommand{\chaptername}{Unit}
\chapter{Introduction}
The objective of this course is to give an introduction to probability theory via the theory of measure and integration. Some amount of inequalities related to probability are also discussed.

This report is divided into five units. The second unit deals with some preliminary concepts related to measure theory, the third unit deals with integration theory where focus is given on Riemann-Steiltjes integrals and Lebesgue integrals. The fourth unit deals with probability spaces and the concept of a random variable as an integral. Finally the fifth unit gives some basic inequalities related to probability.

Almost all the proofs are omitted in this report. Many proofs can be found in either \cite{aps} or \cite{roy}.

\chapter{Preliminaries and Basics of Measure Theory}\label{ch1}

In this Unit we shall look into some preliminary concepts of measure theory, with a few proofs and some examples.

\section{Preliminary Results}

We begin with the following important results which we shall be coming back from time to time.

\begin{property}[Archimedean Property\index{Archimedean Property}]\label{archimedean}
For any pair of positive real numbers $a$ and $b$, there is a natural number $n$ for which $na>b$.
\end{property}

\begin{thm}[Heine-Borel\index{Heine-Borel Theorem}]\label{heine}
Let $F$ be an open covering of a closed and bounded set $A$ in $\mathbb{R}^n$, then a finite subcollection of $F$ also covers $A$.
\end{thm}

Now we have the following

\begin{thm}
An open set can be expressed as union of atmost countable disjoint intervals.
\end{thm}

\begin{proof}
We take $x~\in ~O$ , then there exists an $\epsilon >0$  such that $(x-\epsilon,x+\epsilon)\subset O$.

Now we consider the sets $y_1=\sup\{z>x:(x-\epsilon,z)\subset O\}$ and $y_2=\inf\{z<x:(z,x+\epsilon)\subset O\}$. It is clear that $x\in(y_1,y_2)\subset O$.

For if $x\neq y$, let $O_x$ and $O_y$ be the largest open intervals containing $x$ and $y$ such that $O_x,O_y \subset O$, and either $O_x=O_y$ or $O_x\cap O_y=\phi$.

Now by Propery \ref{archimedean}, there exists an $r_1\in O_x, r_2\in O_y,r_1\neq r_2,r_1,r_2\in Q,$ where $O_x\cap O_y=\phi.$ Since $Q$ is countable, so the no of intervals is countable.

\end{proof}

From the above we immediately have the following

\begin{cor}
\begin{itemize}
\item[(a)] Open sets are Borel sets.
\item[(b)] Closed sets are Borel sets.
\end{itemize}
\end{cor}

It is now clear that $(a,b]$ is a Borel set because

$$(a,b]=\displaystyle\bigcap_{n\in N} (a,b+\frac{1}{n}).$$

Again the set $[a,b)$ is also a Borel set because

$$[a,b)=\cap_{n\in N}(a-\frac{1}{n},b).$$

\section{Length of a Borel set}

It is actually possible to define a notion of length for every Borel set.

\begin{definition}
Let $\mathcal{F}$ be the class of sets which can be expressed as union of finite number of disjoint semi-open intervals.
\end{definition}

As a convention the sets $(a,\infty),(-\infty,b]$ are considered as semi-open intervals.

Cleary the set $\mathcal{F}$ so defined above is a field. We have a natural definition of length of sets in $\mathcal{F}$. If $B\in \mathcal{F}$, then
\begin{align*}
 l(B) &= \text{length of} ~B\\
&=\sum_{i=1}^{n} l(a_{i},b_{i}],
\end{align*}
where  $B=\bigcup_{i=1}^{n} (a_{i},b_{i}] $ with $(a_{i},b_{i}] \cap(a_{i},b_{i}]$.

\begin{definition}
 If $\{(a_{i},b_{i})\}_{i=1,2,...}$ is a sequence of non-overlapping or disjoint semi-open intervals, then
\begin{align*} 
 l(\cup_{i=1}^{\infty}(a_{i},b_{i}])=\sum_{i=1}^{\infty}l(a_{i},b_{i}]
\end{align*}
where $l$ denotes the length.

\end{definition}

If $(a,b]=\cup_{i=1}^{\infty}(a_{i},b_{i}]$, where $\{(a_{i},b_{i})\}_{i=1,2,...}$ is a sequence of disjoint semi-open intervals, then
\begin{align}
 b-a&=\sum_{i=1}^{\infty}l((a_{i},b_{i}])\nonumber\\&=\lim_{n\to\infty}\sum_{i=1}^{n}l((a_{i},b_{i}]).
\end{align}

Given that $\cup_{i=1}^{\infty}(a_{i},b_{i}]=(a,b]$, where $(a_{i},b_{i}]\cap(a_{j},b_{j}]=\phi$.\\ So for any $k$, with $a<a_{k}, b_{k}<b$, we have
\begin{align*}
 b-a>\sum_{i=1}^{k}(b_{i}-a_{i})=\sum_{i=1}^{k}l((a_{i},b_{i}]).
\end{align*}
Taking $k\to\infty$, we have
\begin{align*}
 \sum_{i=1}^{\infty}l((a_{i},b_{i}])\leq b-a.
\end{align*}
Now suppose
\begin{align*}
 (a,b]\subset\cup_{i=1}^{\infty}(a_{i},b_{i}]=\cup_{i=1}^{\infty}I_{i}.
\end{align*}
We now show that $b-a\leq \sum_{i=1}^{\infty}l(I_{i})$.

For any small $\epsilon>0$,
\[
 [a+\epsilon/2,b]\subset(a,b]\subset\cup_{i=1}^{\infty}(a_{i},b_{i}]\subset\cup_{i=1}^{\infty}(a_{i},b_{i}+\epsilon/2^{i+1}).
\]

Then by the Heine-Borel Theorem, there exists a $k$ such that

\begin{align*}
 [a+\epsilon/2,b]&\subset\cup_{i=1}^{k}(a_{i},b_{i}+\epsilon/2^{i+1})\\
\textit{i,e.,}~b-a-\epsilon/2&\leq \sum_{i=1}^{k}(b_{i}-a_{i})+\epsilon\sum_{i=1}^{k}\frac{1}{2^{i+1}}<\sum_{i=1}^{k}(b_{i}-a_{i})+\epsilon/2\\
\textit{i,e.,}~~~~~~b-a<&\sum_{i=1}^{\infty}(b_{i}-a_{i})+\epsilon.
\end{align*}
Since $\epsilon$ is arbitrary, so Equation 2.0.1 is correct.

From the above we get the following as a consequence

\begin{prop}
 If $[a,b]\subset\cup_{i=1}^{k}(a_{i},b_{i})$, then $$b-a\leq \sum_{i=1}^{k}(b_{i}-a_{i}),k\geq 1.$$
\end{prop}

Let $\mathcal{F}$ be the class of all $Borel$ sets such that each set can be expressed as union of atmost countable disjoint semi-open intervals (say ${I_{i}}$) \textit{i.e.,} it is the collection $$\{B~:~B~is~the ~Borel~sets~and~B~is~the~union~of~atmost~countable~semi-open~intervals.\}$$

We have the following basic properties of $l$

\begin{enumerate}
\item $l$ is non-negative.
\item $l$ is finite additive on $\mathcal{F}$, \textit{i.e.,} if $A, B \in \mathcal{F}$, then $l(A\cup B)=l(A)+l(B)$, where $A\cap B=\phi.$ 
\item $l$ is countably addtive on $\mathcal{F}$\textit{i.e.,} if $\{B_{1},B_{2},\dots\}$ is a sequence of sets in $\mathcal{F}$ such that 
 $B_{i}\cap B_{j}=\phi$,then $l(\cup_{i=1}^{\infty}B_{i})=\sum_{i=1}^{\infty}l(B_{i})$.
\end{enumerate}

\begin{definition}[Borel $\sigma$-field\index{Borel $\sigma$-field}]
The smallest $\sigma$-field containing all open intervals and semi open intervals is called the Borel $\sigma$-field.
\end{definition}

\begin{thm}[Caratheodory Extension Theorem\index{Caratheodory Extension Theorem}]
Let $\mathcal{F}$ be a class (say a field) containing all semi-open intervals such that a function 
$\mu$ on $\mathcal{F}$ satisfies the non-negativity and countable subadditivity properties. Then $\mu$ can be extended to a measure on the Borel $\sigma$-field. This extension is unique.
\end{thm}

The above theorem ensures that $l$ extends to a measure $\lambda$ on the $Borel ~\sigma$-field, where $\lambda$ is called the Lebesgue measure\index{Lebesgue measure} which has the following properties.

\begin{itemize}
\item[1.] Lebesgue measure of a finite set is zero.
\item[2.] $\lambda(x+cA)=c\lambda(A)$.
\end{itemize}

It is important to note the following

\begin{enumerate}
\item Probability is a finite measure.
\item Lebesgue measure is not finite, but it is $\sigma$-finite as $\mathbb{R}=\cup_{n\in I}(n,n+1]$.
\end{enumerate}

\begin{definition}
Let $\mathcal{C}$ be a class of subsets of $\mathbb{R}$. Then $\sigma(\mathcal{C})$ is the smallest $\sigma$-field 
containing the sets in $\mathcal{C}$, or $\sigma(\mathcal{C})$ is generated by $\mathcal{C}$.
\end{definition}
\begin{thm}
 For any class of sets $\mathcal{C}, \mathcal{C}\cap B=\{A\cap B: A\in \mathcal{C}$\}. Let B be any Borel set. Then $\sigma(\mathcal{C})\cap B=\sigma(\mathcal{C}\cap B)$.
\end{thm}
\begin{definition}[Counting Measure\index{counting measure}]
Let S be a finite or countable set and $\mathcal{F}$ be a field or $\sigma$-finite of subsets on S. Then $\mu$ on (S, $\mathcal{F}$) to be\\ 
$\mu(B)=\left\{\begin{array}{cc}
              \text{number of elements in $B$}, & \text{for any finite set $B \in \mathcal{F}$} \\
\infty, & \text{otherwise}
             \end{array}
\right.$
\end{definition}

\begin{rem}
Counting measure and Lebesgue measure are mutually sigular and length of an unbounded set is infinite.
\end{rem}

\section{Distribution Functions}

We now generalize the notion of length in the following manner

\begin{definition}
Let $F:\mathbb{R}\longrightarrow (0,\infty)$ be a right continuous monotonically non-decreasing function $F$, we define a notion of length as follows:
\begin{align*}
 l(a,b]&=F(b)-F(a)\\
l(a,b)&=F(b-)F(a)\\
l[a,b)&=F(b-)-F(a-)\\
l[a,b]&=F(b)-F(a-)
\end{align*}
where $F(x-)=\displaystyle\lim_{\begin{array}{cc}
  y\to x & \\
y<x & 
 \end{array}
} F(y)$
\end{definition}

\begin{thm}
 If $(a,b]\subset\cup_{i=1}^{\infty}I_{i}$, where $I_{i}=(a_{i},b_{i}]$ and $I_{i}\cap I_{j}=\phi$, then $$F(b)-F(a)\leq\sum_{i=1}^{\infty}\{F(b_{i})-F(a_{i})\}.$$
\end{thm}

\begin{proof}
Let $\epsilon>0$ be given. Then we have $$[a+\delta,b]\subset(a,b]\subset\cup_{i=1}^{\infty}(a_{i},b_{i}]\subset\cup_{i=1}^{\infty}(a_{i},b_{i}+\delta_{i})$$ where $F(a+\delta)-F(a)\leq\epsilon/2$ and $F(b+\delta_{i})-F(b)\leq\epsilon/2^{i+1}$.

By the Heine-Borel Theorem, there exists an $N$ such that \begin{align*}
                                              [a+\delta,b]&\subset\cup_{i=1}^{N}(a_{i},b_{i}+\delta_{i})\\
\Rightarrow F(b)-F(a+\delta)&\le\sum_{i=1}^{N}(F(b_{i}+\delta_{i})-F(a_{i}))\\
\Rightarrow F(b)-F(a+\delta)&\le\sum_{i=1}^{N}(F(b_{i}+\delta_{i})-F(b_{i}))+\sum_{i=1}^{N}(F(b_{i})-F(a_{i}))\\
&\le\sum_{i=1}^{N}\epsilon/2^{i+1}+\sum_{i=1}^{N}(F(b_{i})-F(a_{i}))\\
&\le\epsilon/2+\sum_{i=1}^{N}(F(b_{i})-F(a_{i})).
                                             \end{align*}
Now \begin{align*}
     F(b)-F(a)&=F(b)-F(a+\delta)+F(a+\delta)-F(a)\\
&\le F(b)-F(a+\delta)+\epsilon/2\\
&\le \sum_{i=1}^{\infty}(F(b_{i})-F(a_{i}))+\epsilon ~~\forall \epsilon>0.
    \end{align*}
Since $\epsilon$ is arbitrary, hence it completes the rest of the proof.
\end{proof}

\begin{lem}
If $[x , y]\subset\cup_{i=1}^{N}(a_{i},b_{i})$, then $F(y)-F(x)\le\sum_{i=1}^{N}(F(b_{i})-F(a_{i})),$ where F is a non-decreasing function. 
\end{lem}

\begin{definition}[Probability Distribution\index{probability distribution}]
A probability distribution is a (probability) measure on $(\mathbb{R},\mathcal{B})$. 
\end{definition}

\begin{note}
$F_0$ = the field of all possible finite disjoint union of semi-open (right closed) intervals.\\
$B\in F_0 \Rightarrow B=\cup_{i=1}^{k}(a_{i},b_{i}], k\ge 1$, with $(a_{i},b_{i}]\cap(a_{j},b_{j}]=\phi.$
\end{note}

\begin{definition}
 The definition of length on $F_0$ to be\begin{enumerate}
                                         \item $l(B)=\sum_{i=1}^{k}(F(b_{i})-F(a_{i}))$, where $B=\cup_{i=1}^{k}(a_{i},b_{i}]$,
\item $l(B)=\sum_{i=1}^{\infty}(F(b_{i})-F(a_{i}))$, where $B=\cup_{i=1}^{\infty}(a_{i},b_{i}]$.
                                        \end{enumerate}
\end{definition}

The above definition holds due to the following theorem.

\begin{thm}
 If $(a,b]=\cup_{i=1}^{\infty}(a_{i},b_{i}]$ such that $(a_{i},b_{i}]\cap(a_{j},b_{j}]=\phi$, then $F(b)-F(a)=\sum_{i=1}^{\infty}(F(b_{i})-F(a_{i}))$.
\end{thm}

\begin{definition}[Distribution Function\index{distribution function}]
A distribution function $F$ is non-negative monotonically non-decreasing right continuous function such that $F(-\infty)=0, F(\infty)=1$, where $F:\mathbb{R}\longrightarrow [0,1]$.
\end{definition}

Given a distribution function $F$ we have 

\begin{align*}
 l({x})&=l(\cap_{n=1}^{\infty}(x-\frac{1}{n},x])\\
&=\lim_{n\to\infty}l{(x-\frac{1}{n},x]}\\
&=\lim_{n\to\infty}{F(x)-F(x-\frac{1}{n},x)}\\
&=F(x)-F(x-)\\
&=\text{Jump or probability mass at x.}
\end{align*}

\begin{definition}[Continous distribution\index{continous distribution}]
If the distribution function is continuous, then the probability distribution is called the continuous probability distribution.
\end{definition}

\begin{note}
For any continuous distribution, probability or measure of any interval $(a,b)$ or $(a,b]$ or $[a,b]$ or $[a,b)$ is always $F(b)-F(a)$.
\end{note}

\begin{thm}
The number of discontinuity points of a monotonic function is atmost countable.
\end{thm}

\begin{note}
For any distribution function $$0 \le F(x)-F(x-)\le F(x)\le F(\infty)=1.$$
\end{note}

\begin{definition}[Purely discrete distribution\index{purely discrete distribution}]
Let $S$ be the set of discontinuity points of a distribution function $F$ with $S$ is atmost countable. If $\sum_{x\in S}{F(x)-F(x-)}=1$, the probability distribution with distribution $F$ is said to 
be {\bf purely discrete distribution} or {\bf atomic distribution}\index{atomic distribution} where $S$ is called the support of the distribution and its elements are called the {\bf mass points}\index{mass points}. 
\end{definition}

\begin{rem}
We denote a probability measure by $P$ and uniquely determined by the distribution function $F$. Also we have $P\{(a,b]\}=F(b)-F(a)$.
\end{rem}

\begin{lem}
\begin{enumerate}
\item If $S$ is the support of a probability distribution $P$, then $P(\{y\})=0$, for $y\notin S, i.e., P(S^c)=0$.
\item If $P(\{y\})>0$, then $P(S\cup \{y\})>1$.
\item If $P$ is a purely discrete distribution with support $S$, then\\$P(B)=P(B\cap S)=\sum_{x\in S\cap B}P(\{x\})=\sum_{x\in S\cap B}\{F(x)-F(x-)\}$.
\end{enumerate}
\end{lem}

\begin{prop}[Regularity Property\index{regularity property}]
Let $P$ be probability distribution on $(\mathbb{R},\mathcal{B})$. Then for any Borel set B and $\epsilon>0,~ \exists$ a closed set C and an open set O such that 
$C\subset B\subset O$ and $P(O-C)<\epsilon$.
\end{prop}

If $P$ is a probability distribution such that $P(B)=0$ or $P(B)=1$ for every Borel set $B$, then $P(\{x_0\})=1$. For every probability distribution $P$ on $(\mathbb{R},\mathcal{B})$, there exists a smallest closed set $S$ such that $P(S)=1$.

\begin{rem}
Such an $S$ is called the support\index{support} of the distribution.
\end{rem}

\chapter{Integration Theory}

Integration theory is a very important aspect of all of analysis sepcially in probabilty theory as we shall see in this unit. The normal school level integration that we do is what is called \textit{Riemann integration}\index{Riemann integration}, in this section we shall get accquanted with other theories of integration which shall be very useful in our study of probability theory.

\section{Riemann Integration}

In our usual Riemann integration the necessary condition to define $\int_{a}^{b} f(x) dx$, is that the function $f$ should be bounded while a sufficient condition to define $\int_{a}^{b} f(x) dx$, is that the function $f$ should be continuous.

\begin{rem}
Boundedness is not sufficient condition to define Riemann Integration. Here is the counter example of a function which is bounded but not Riemann Integrable,\\
$f(x)=\left\{ \begin{array}{cc}
               1, & \text{x is rational}\\
0, & \text{x is irrational.}
            \end{array}
\right.$
\end{rem}

\begin{note}
Let $P=\{x_0,x_1,\dots,x_n\}$, where $x_0=a, x_n=b$ and $x_0<x_1<\dots<x_n$. Let $P_1,P_2$ be two partition of $[a,b]$. We say $P_2$ is a refinement of $P_1~ if ~P_2\supset P_1$.
\end{note}

\begin{definition}[Upper Sums, Lowe Sums\index{upper sum}\index{lower sum}]
Suppose $f$ is bounded on $[a,b]$, that is there exists an $m, M$ such that $m<f(x)<M, \forall x\in[a,b]$. Given a partition $P=\{x_0=a,\dots\ x_i,\dots,x_n=b\}$. Let\begin{align}
    m_i&=\inf_{x\in[x_{i-1},x_{i}]} f(x)\label{d1}\\ 
M_i&=\sup_{x\in[x_{i-1},x_{i}]}f(x)\label{d2}.
   \end{align}
Then\begin{align*}
     L(P)&=\sum_i (x_{i}-x_{i-1}) m_i,\\ L(P)&=\sum_i (x_{i}-x_{i-1}) M_i.    
    \end{align*}
Also, for each $P$, we have \begin{align*}
                             L(P)\le U(P).
                            \end{align*}
\end{definition}

The following are some of the important properties of the upper sums and lower sums that we have just defined.

\begin{lem}
Let $P_1, P_2$ be any two partition on $[a,b]$. Then $L(P)\le U(P)$.
\end{lem}

\begin{lem}
If $P_1\subset P_2$, then $L(P_1)\le L(P_2)~and~U(P_1)\le U(P_2)$.
\end{lem}

\begin{proposition}
If a function is continuous, then that function is Reimann Integrable.
\end{proposition}

\begin{thm}
 If $f$ is continuous on $[a,b]$, then $f$ is uniformly continuous on $[a,b]$.
\end{thm}

\begin{prop}
If $f$ is bounded on $P$, then \begin{align*} \sup_{x\in P}f(x)-\inf_{x\in P}f(x)\ge\sup_{x,y\in P}|f(x)-f(y)| \end{align*}
\end{prop}

\begin{note}
We denote the norm of a partition\index{partition norm} as the length of the greatest subinterval.
\end{note}

\begin{thm}
A Reimann Integrable function can have atmost countably many discontinuity points.
\end{thm}

\begin{exm}
$f(x)=\left\{\begin{array}{cc}
              0, & 0<x<1/2\\
1, & 1/2\le x<1
            \end{array}
\right.$\\
We have $\int_{0}^{1} f(x)dx=\int_{1/2}^{1}dx=\frac{1}{2}$. Therefore the function is Reimann Integrable. 
\end{exm}

We shall now state some important properties of integration.

\begin{enumerate}
                                 \item $\int_{a}^{b} f(x) dx$ is linear, that is 
                                 \begin{align*}
                                                                                 \int_{a}^{b}\{c f(x)+d g(x)\} dx = c \int_{a}^{b} f(x)dx + d \int_{a}^{b} g(x) dx.
                                                                                \end{align*}
				 \item If $f(x)\le g(x)  ~\forall ~x\in [a,b]$, then\begin{align*} 
				                                                   \int_{a}^{b} f(x) dx\le \int_{a}^{b} g(x) dx.
										     \end{align*}

                                \end{enumerate}
                                
\section{Riemann-Steiltjes Integration}

After the preliminaries on Riemann integration, we now move on to another theory of integration called the Riemann-Stieltjes integration theory. We begin with the following definitions of upper and lower sums in the context of Riemann-Steiltjes integration. \index{upper sum} \index{lower sum}

\begin{definition}
 To define $\int_{a}^{b} f(t)dF(t)$, we assume $(i)~f$ is bounded, $(ii)~F$ is monotonically non-decreasing. For a partition, $P=\{x_0=a,x_1,\dots,x_n=b\}$, \label{P:a}  we have  $F(x_i)-F(x_{i-1})\ge 0.$ And we define Upper sum\index{upper sum} and Lower sum\index{lower sum} as\begin{align*}
                            U(P,f,F)=\sum_{i=1}^{n}M_i \{F(x_i)-F(x_{i-1})\},\\ \label{a}
L(P,f,F)=\sum_{i=1}^{n}m_i \{F(x_i)-F(x_{i-1})\}. 
                           \end{align*}
where $m_i$ and $M_i$ are defined by (\ref{d1})   and (\ref{d2}) respectively. For any partition $P_1, P_2$, we have \begin{align*}
                                       L(P,f,F)\le U(P,f,F).
                                      \end{align*}
Moreover, if $P_1\subseteq P_2$, then \begin{align*}
                                     U(P,f,F)&\ge U(P,f,F)~\text{and} \\
L(P,f,F) &\le L(P,f,F).
                                    \end{align*}
\end{definition}

We are now in a position to define what is a Riemann-Steiltjes integral.

\begin{definition}[Riemann-Steiltjes Integral\index{Riemann-Steiltjes integration}]
For every $\epsilon >0$, there exists a partition $P_\epsilon$ of $[a,b]$ such that for every refinement $P$ of $P_\epsilon, ~|U(P,f,F)-L(P,f,F)|<\epsilon$. In that case $f\in R(F)[a,b]$ that is $f$ is Reimann-Stieltjes Integrable with respect to $f$ on $[a,b].$
\end{definition}

We also have the following alternate definiton of the above

\begin{definition}
Suppose $f\in R(F)[a,b]$, then \begin{align*}
                                                                                                         \lim_{\delta\to 0} S(P_\delta,f,F)&=\int_{a}^{b} f(x)dF(x),
                                                                                                        \end{align*}
where \begin{align*}
      ~~~~~~~~~ S(P_\delta,f,F)&=\sum_{i=1}^{k}f(t_i)\{F(x_i)-F(x_{i-1})\}
      \end{align*}
where $t_i$ is a point in the $i^{th}$ subinterval of the partition $P$ of Definition \ref{P:a}.
\end{definition}

\begin{note}
For any partition P, we have $L(P,f,F)\le S(P,f,F)\le U(P,f,F).$
\end{note}

We now state the following properties of Riemann-Steiltjes integration.

\begin{enumerate}
                                                   \item If $f,g\in R(F)[a,b],$ then $cf+g\in R(F)[a,b].$
\item If $f\le g$ on $[a,b]$, then \begin{align*}
                                    \int_{a}^{b} f(x)dF(x)\le \int_{a}^{b}g(x)dF(x).
                                   \end{align*}
\item Suppose $F,G$ are two monotonic functions, if $f\in R(F)[a,b]$ and also $f\in R(G)[a,b]$ then $$\int_{a}^{b}f(x)d(F+G)(x)=\int_{a}^{b}f(x)dF(x)+\int_{a}^{b}f(x)dG(x).$$

\item $\int_{a}^{b}f(x)dF(x)=-\int_{a}^{b}f(x)dF(x)$.
\end{enumerate}

The following theorem is one of the most important results in integration theory.

\begin{thm}[Integration by parts\index{integration by parts}]
If $f\in R(F)$ on $[a,b]$, then $F\in R(f)$ on $[a,b]$ and $$\int_{a}^{b} f(x)dF(x)+\int_{a}^{b} F(x)df(x)=F(b)f(b)-F(a)f(a).$$
\end{thm}

\begin{thm}
If $f$ is continuous on $[a,b]$ and $F$ is monotonically non-decreasing, then $f\in R(F)$ on $[a,b]$.
\end{thm}

\begin{ex}
 If $f$ is a continuous function on $[a,b]$ and $F(x)=\left\{\begin{array}{cc}
              0, & 0\le x<\frac{1}{2},\\
1, & \frac{1}{2}\le x\le 1.
            \end{array}
\right.$\\
Then find the value of  $\int_{0}^{1}f(x)dF(x)=\lim_{\delta\to 0}S(P,f,F)$.
\end{ex}

\begin{lem}
Let $g$ be a continuous strictly monotonic function on $[a,b]$, then $g\{[a,b]\}=\{g(y):a\le y\le b\}$ is an interval.
\end{lem}

\begin{lem}
A strictly monotonic continuous function $g ~on ~[a,b]$ is invertible.
\end{lem}

\begin{lem}
The function $g^{-1}$ is also a continuous function.
\end{lem}

The above three lemmas give us the following important result.

\begin{thm}[Change of variables\index{change of variables}]
Let $f\in R(F)$ on $[a,b]$ and let $g$ be a strictly monotonic continuous function defined on an interval $S$ having endpoints $c$ and $d$. Assume that $a=g(c),~b=g(d).$ Let $h$ and $G$ be the 
composite function defined as follows\begin{align*}
                                       h(x)=f[g(x)], ~~G(x)=F[g(x)],~~~~ if x\in S.
                                      \end{align*}
Then $h\in R(G)$ on $S$ and we have $\int_{a}^{b}f dF=\int_{a}^{b}h dG$. That is,\begin{align*}
                                                                                  \int_{g(c)}^{g(d)}f(t)dF(t)=\int_{c}^{d}f[g(x)]d\{F[g(x)]\}.
                                                                                 \end{align*}
\end{thm}

\begin{thm}
 Given $a<c<b$. Define F on $[a,b]$ as $F(x)=\left\{\begin{array}{cc}
              F(a), & a\le x<c\\
F(b), & c< x\le b
            \end{array}
\right.$. Suppose $f$ is a function on $[a,b]$ such that atleast one of the function on $F$ is a right continuous at $c$ and the other is left continuous at c, then $f\in R(F)$ on $[a,b]$ and $$\int_{a}^{b} f(x)dF(x)=f(b)\{F(c+)-F(c-)\}.$$
\end{thm}

\begin{thm}
Any finite sum $\sum_{k=1}^{n}a_k$ can be written as \begin{align*}
                                                       \int_{0}^{n}f(x) d[x]
                                                      \end{align*}
where $f(x)=a_k,~k-1<x\le k, ~k=1,2,\dots,n ~;f(0)=0~and~[x]$ is the greatest integer less than or equal to x.
\end{thm}

\begin{thm}[Euler Summation Formula\index{Euler summation formula}]
If $f$ has a continuous derivative $f'$ an $[a,b]$, then we have\begin{align*}
                                                                                             \sum_{a<n\le b}f(n)=\int_{a}^{b}f(x)dx+\int_{a}^{b}f'(x)(x-[x])dx+f(a)(a-[a])-f(b)(b-[b]).
                                                                                            \end{align*}
\end{thm}

We can prove the above result using the change of variable theorem and the two other theorems that follow it.

\begin{thm}
Assume $f\in R(F)$ on $[a,b]$ and assume that $F$ has a continuous derivative $F'$ on $[a,b]$. Then the Riemann integral $\int_{a}^{b}f(x)F'(x)dx$ exists and we have\begin{align*}
  \int_{a}^{b}f(x) dF(x)=\int_{a}^{b}f(x)F'(x)dx.
 \end{align*}
\end{thm}

\section{Lebesgue Integration}

Having learnt Riemann and Riemann-Steiltjes integration we shall now go into another very important theory of integration called Lebesgue integration\index{Lebesgue integration}, which was devised by Henri Lebesgue, a French mathematician in the twentieth century.

To motivate the need for this new theory of integration we look at the following function called Dirichlet's function\index{Dirichlet function}\\

$f(x)=\left\{ \begin{array}{cc}
               1, & \text{x is rational,}\\
0, & \text{x is irrational.}
            \end{array}
\right.$\\

We know that for any monotonic function $F$, $\int_{0}^{1}f(x)dF(x)$ is not defined in Riemann-Stieltjes sense.

Also we have $f(x)=1\times I_Q(x)+0\times I_{Q^c}(x)$ where $I$ stands for indication or characteristic function\index{characteristic function}\index{indication function}, defined by $I_Q(x)=\left\{ \begin{array}{cc}
               1, & x\in Q,\\
0, & x\notin Q.
            \end{array}
\right.$

\begin{definition}[Simple function\index{simple function}]
 A function $f$ is said to be simple if it is measurable and its range is a finite set.
\end{definition}

\begin{definition}
Let $(\mathbb{R},\mathcal{B})$ be a measurable space, where $\mathcal{B}$ is a $Borel~\sigma$-field. Let $\lambda$ be a measure on $(\mathbb{R},\mathcal{B})$. Consider the map $f:\mathbb{R}\longrightarrow E$, where
$E$ is a finite subset of $\mathbb{R}$ and $f^{-1}(B)=\{x\in \mathbb{R}:f(x)\in B\}\subset \mathcal{B}$. Then we have \begin{align*}
                                                                                                                       \int_{B}f d\lambda&=\int_{B}(1\times I_{Q}(x)+0\times I_{Q^c}(x))d\lambda(x)\\
&=1\times\lambda(B\cap Q)+0\times\lambda(B\cap Q^c).
                                                                                                                      \end{align*}
\end{definition}

\begin{definition}
If $f$ is simple function with range $\{a_1,\dots,a_k\}$ on a measurable space $(\mathbb{R},\mathcal{B},\lambda)$, then $$\int_B fd\lambda=\sum_{i=1}^{k}a_i\lambda(B\cap A_i)$$
where $A_i=\{x\in\mathbb{R}:f(x)=a_i\}$.
\end{definition}

\begin{definition}
We have \begin{align*}f(x)&=f^+(x)-f^-(x)\\
          f^+(x)&=max\{f(x),0\}\\
f^-(x)&=max\{-f(x),0\}\\
\int_Bfd\lambda&=\int_Bf^+d\lambda-\int_Bf^-d\lambda.
         \end{align*}
\end{definition}

\begin{rem}
In general $\lim_{n\to\infty}\int_{0}^{1}f_n(x)dx=\int_{0}^{1}f(x)dx$ is not true. Following example shows that this condition is violated.
\end{rem}

\begin{exm}
  $f_n(x)=\left\{\begin{array}{cc}
              n, & 0\le x<\frac{1}{n}\\
1, & \frac{1}{n}< x\le 1
            \end{array}
\right.$.
In this example $\lim_{n\to\infty}\int_{0}^{1}f_n(x)dx\ne\int_{0}^{1}f(x)dx$.
\end{exm}

\begin{rem}
If $$\mid f_n(x)\mid \le g(x) ~~(free~of~n)$$
and $\int_{0}^{1}g(x)dx$ is finite, then $$\lim_{n\to\infty}\int_{0}^{1}f_n(x)dx=\int_{0}^{1}f(x)dx.$$
\end{rem}

We now state a very important result in the theory of Lebesgue integration.

\begin{thm}[Dominated Convergence Theorem\index{Dominated convergence theorem}]
Suppose $\{f_n\}$ is a sequence of integrable functions on $(\mathbb{R},\mathcal{B},\lambda)$ such that $f_n(x)\to f(x)$ for almost all x 
and $\mid f_n(x)\mid \le g(x)$ for almost all x such that $g$ is an 
integrable function. Then $f$ is integrable and $\int fd\lambda=\lim_{n\to\infty}\int f_nd\lambda$.
\end{thm}

The following example illustrates the need for the boundedness condition in the Dominated Convergence Theorem.

\begin{exm}
Define $f:[0,\infty]\longrightarrow [0,\infty]$ by \\
 $f(x)=\left\{\begin{array}{cc}
              0, & 0\le x\le 1\\
1, & 1< x\le 1+\frac{1}{2^3}\\
0, & 1+\frac{1}{2^3} <x\le 2\\
2, & 2<x\le 2+\frac{1}{2^3}\\
0, & 2+\frac{1}{2^3}<x\le 3\\
\vdots\\
n, & n<x\le n+\frac{1}{(n+1)^3}\\
\vdots
            \end{array}
\right.$\\
Here we have \begin{align*} 
     \int_{0}^{\infty}f(x)dx=\sum_{n=1}^{\infty}n.\frac{1}{(n+1)^3}\le \sum_{n=1}^{\infty}\frac{1}{(n+1)^2}<\infty.
    \end{align*}
The above series is convergent, but the function is unbounded.
\end{exm}

Below we state some important results without proof. Almost all the results below follow from definitions or from the results stated earlier to them.

\begin{thm}
 Suppose $f$ is a non-negative finite valued measurable function on a measure space $(\Omega,\mathcal{F},\mu)$. Then there exists a sequence on non-negative simple functions $\{f_n\}_{n=1,2,\dots}$
such that $f_n\le f_{n+1}$ (pointwise) and if $f(x)<\infty ~\forall~x\in\Omega$, then $f_n(x)\to f(x)~~\forall ~x.$
\end{thm}

\begin{lem}
 If $f$ and $g$ are two simple functions on $(\Omega,\mathcal{F},\mu)$ such that $f(x)\le g(x)$. Then $\int fd\mu\le\int gd\mu.$
\end{lem}

\begin{rem}
 For any non-negative measurable function $f$ on $(\Omega,\mathcal{F},\mu)$ there exists a sequence $\{f_n\}_{n=1,2,\dots}$ of non-negative simple functions such that $f_n\le f_{n+1}$ and 
$f_n\to f$ a.e. Then $\int f_n d\mu \le \int f_{n+1}d\mu$.
\end{rem}

\begin{definition}
For any non-negative measurable function $f$ on $(\Omega,\mathcal{F},\mu)$ satisfying $\mu(\{x:f(x)=\infty\})=0$, define $\int_B fd\mu=\lim_{n\to\infty}\int_B f_n d\mu$ for $B\in \mathcal{F}$
\end{definition}

\begin{definition}
 For any measurable function $f$ on $(\Omega,\mathcal{F},\mu)$ satisfying $\mu(\{x:f(x)=\pm\infty\})=0$, define $\int_B fd\mu=\int_B f^+ d\mu-\int_B f^-d\mu$ provided $\int_\Omega f^+d\mu$ or $\int_\Omega f^-d\mu<\infty$ where $f^+$ and $f^-$ are integrable functions.
\end{definition}

Below we list some of the important properties of Lebesgue integrals.

If $f,g$ are integrable functions on $(\Omega,\mathcal{F},\mu)$, then
\begin{enumerate}
                  \item $\int(af+g)d\mu=a\int fd\mu+\int gd\mu$.
\item $|\int fd\mu|\le \int|f|d\mu$.
\item If $f\le g$, then $\int fd\mu\le\int gd\mu$.
\item If $\int_B fd\mu\le\int_B gd\mu~~\forall~B\in\mathcal{F}$, then $f\le g$ almost everywhere.
\item If $\int_B fd\mu=\int_Bd\mu~~\forall~B\in\mathcal{F}$, then $f=g$ almost everywhere.
                 \end{enumerate}
The following are some easy consequences of the above properties.
         
\begin{lem}
If $f$ be a simple function on $(\Omega,\mathcal{F},\mu)$ such that $\int_B fd\mu\ge0~~\forall~B\in\mathcal{F}$, then $f\ge0$.
\end{lem}

\begin{cor}
If $f$ is simple function on $(\Omega,\mathcal{F},\mu)$ such that $\int_B fd\mu=0~~\forall~B\in\mathcal{F}$, then $f=0$.
\end{cor}

\begin{lem}
If $f$ is an integrable function on $(\Omega,\mathcal{F},\mu)$, then $\exists$ a sequence of simple functions $\{f_n\}_{n=1,2,\dots}$ such that $\int_B f d\mu=\lim_{n\to\infty}\int_B f_n d\mu$.
\end{lem}

\begin{definition}\index{integrability}
A measurable function $f$ on $(\Omega,\mathcal{F},\mu)$
 is said to be integrable if $$\int f_+ d \mu < \infty~~\text{and}~~\int f_- d \mu < \infty$$
$$\text{or},~~~ \int |f|d\mu <\infty \Leftrightarrow f \in L_1(\Omega,\mathcal{F},\mu) $$
$$\text{or},~~~ \int fd\mu ~~\text{exists}~~~\Leftrightarrow \int f_+ d \mu < \infty ~~\text{or}~~\int f_- d \mu < \infty. $$
\end{definition} 
\begin{proposition}
 If $f$ is integrable, then $f$ is finite almost everywhere, i. e. $\mu \big\{ \{x\} : |f(x)|=\infty \big\}=0.$
\end{proposition}

\begin{proposition}
 If $f \geq 0$, then $\int fd\mu =0$ implies $f=0$  almost everywhere.
\end{proposition}

\begin{thm}
Suppose that $f$ and $g$ are two integrable  functions on $(\Omega,\mathcal{F},\mu)$ such that $\int_{B} f d\mu \leq \int_{B} g d\mu,~~\forall ~~B\in \mathcal{F} \Rightarrow f \leq g $ almost everywhere.
\end{thm}

\begin{thm}
Suppose that $f$ and $g$ are two integrable  functions on $(\Omega,\mathcal{F},\mu)$ such that $\int_{B} f d\mu = \int_{B} g d\mu,~~\forall ~~B\in \mathcal{F} \Rightarrow f = g $ almost everywhere.
\end{thm}

\begin{thm}
Let $f$ be a nonnegative integrable function and $\int _{\Omega} f d\mu < \infty.$ Define $\lambda (B)=\int _{B} f d\mu.$ Then $\lambda$ is a measure on $(\Omega,\mathcal{F}).$
\end{thm}

\begin{thm}
 Let $f$ be a nonnegative measurable function on $(\Omega,\mathcal{F},\mu).$ Let $\{f_n\}_{n=1}^\infty$ be a sequence of non-negative simple function on $(\Omega,\mathcal{F},\mu)$ such that $f_n(x) \uparrow f(x). $ Then $\lim_{n \to \infty}\int_{\Omega} f_n d\mu=\sup \{ \int_{\Omega} s d \mu: 0 \leq s \leq f \},$ where $s$ is a non-negative simple valued functions on $(\Omega,\mathcal{F},\mu)$.
\end{thm}

\begin{definition}
If $f$ is a non-negative measurable function, then we define
$$\int_{\Omega}f d\mu=\sup\left\{\int_{\Omega} s d \mu: 0 \leq s \leq f \right\},$$ where $s$ is a non-negative simple valued functions on $(\Omega,\mathcal{F},\mu)$.
\end{definition}

\begin{thm}
 If $f \leq g$, then \begin{equation}
                      \int f d\mu \leq \int g d\mu \label{leq}
                     \end{equation}
in the sense that if $\int g d\mu < \infty,$ then $\int f d\mu$ exists and (\ref{leq}) holds. Similarly, if $\int f d\mu > - \infty,$ then $\int g d\mu $ exists and (\ref{leq}) holds.
\end{thm}

We end this discussion with the following famous theorem.

\begin{thm}[Monotone Convergence Theorem\index{monotone convergence theorem}]
Let $\{f_n\}_{n=1}^\infty$ be a sequence of non-negative measurable functions $(\Omega,\mathcal{F},\mu)$ such that 
 $f_n(x) \uparrow f(x)~~\forall x. $ Then $\int f d\mu =\lim_{n \to \infty} \int f_n d\mu.$
\end{thm}

\chapter{Probability Measure and Random Variables}

In this unit we shall study what are probability measures and random variables. The objective here is to demonstrate them as integrals.

\section{Probability Measure}

Let $f$ be a measurable function on $(\Omega,\mathcal{F},\mu)$, where $f : \Omega \to \mathcal{R}.$ The objective is to define a measure on $(\mathcal{R},\mathcal{B})$. If $f$ is a measurable function on $(\mathcal{R},\mathcal{B})$, then $B \in \mathcal{B} \Rightarrow f^{-1}(B) \in \mathcal{F}.$ Thus for $B \in \mathcal{B} $ we define $\mu_0(B)= \mu(f^{-1}(B)).$ Then $\mu_0$ is a measure on $(\mathcal{R},\mathcal{B})$ is called probability measure\index{probability measure}.

\begin{definition}
 If $\mu$ is a probability measure and $f$ is a random variable, then $\mu_0$ is called the probability distribution\index{probability distribution} of $f$, where $\mu_0$ is defined as $\mu_0(B)= \mu(f^{-1}(B)).$
\end{definition}

With the above two definitions and the results from the previous unit we have the following.

\begin{thm}
 Let $f$ be a measurable function on $(\Omega,\mathcal{F},\mu)$ and $g$ be an integrable function on
$(\mathbb{R},\mathcal{B},\mu_0).$ Then $g(f)$ is integrable on $(\Omega,\mathcal{F},\mu)$ and $$\int_{f^{-1}(B)}g(f) d \mu =\int_B g d\mu,$$ where $\mu_0(B)= \mu(f^{-1}(B)).$
\end{thm}

\section{Random Variable and Mathematical Expectation}

\begin{definition}[Random Variable, Mathematical Expectation]
 Let $(\Omega,\mathcal{F},P)$ be a probability space. A random variable\index{random variable} on $X$ on $(\Omega,\mathcal{F},P)$ is a real valued measurable mapping on $(\Omega,\mathcal{F},P).$ Then the mathematical expectation $E(X)$\index{expectation} is defiend as 
\begin{align*}
 E(X)&= \int_\Omega X dP\\
&= \int_{-\infty}^\infty x dF(x)\\
&= \int_{\mathbb{R}}x d \mu,
\end{align*} where $\mu$ is a measure  on $(\mathbb{R}, \mathcal{B})$ defined by the distribution function.
\end{definition}
\begin{rem}
 Given a random variable $X$, we have the distribution function $F$ defined as $F(x)=P(X \leq x)$. Here $g(x)=x$ and $f=X$. So,
$$\int_\Omega X dP= \int_\Omega g(f) dP=\int_{\mathbb{R}}x d\mu.$$
\end{rem}

We now state some important properties of mathematical expectation which can be proved very easily by applying the definition.

 \begin{enumerate}
  \item $E(aX+b)=aE(X)+b.$
\item  \begin{align*}
                X \geq 0,~~   E(X)& \Rightarrow X=0~~\textit{a. e.}\\
\Leftrightarrow P(X>0)&=0\\
\Leftrightarrow P(X=0)&=1.
                  \end{align*}
\item $E(XI_A)=0,~~\forall A \in \mathcal{F}\Rightarrow P(X=0)=1.$
\item $X \leq Y \Rightarrow E(X) \leq E(Y). $
 \end{enumerate}

It is important to note that $X$ is a random variable on a probability space $(\Omega,\mathcal{F},P)$ if and only if $X$ is a measurable function on the measure space $(\Omega,\mathcal{F},P)$.

We also have the following definition.

\begin{definition}
 A random variable $X$ is said to be integrable if $E(|X|)=\int_{-\infty}^\infty | x| dP <\infty $
\end{definition}

\section{Expectations of Some Distributions}
 
We now mention below the mathematical expectations of some well known distributions.

\begin{enumerate}
 \item If $X \sim $ Bernouli distribution\index{Bernoulli distribution} $(P)$, then the distribution function $F$ is defined as
$$F(x)=\left\{\begin{array}{cc}
              0,& x<0,\\
p, & 0\leq x \leq 1,\\
1, & x\geq 1.
             \end{array}
\right.$$
Then $E(X)=\int_{-\infty}^\infty x dF(x)=1-p.$
\item If $X$ is a random
 variable of a purely discrete distribution\index{purely discrete distribution} with support $\{x_1,x_2,\dots,x_n\}$, then $E(X)=\int_{-\infty}^\infty x dF(x)=\sum_{n=1}^kx_i(F(x_i)-F(x_{i-1}))=\sum_{n=1}^kx_iP(X=x_i).$
\item If $X\sim$ Binomial distribution\index{Binomial distribution} $(n,P),$  that means
$P(X=k)={n \choose k}P^k (1-P)^{n-k},~~k=0,1,2,\dots,n,$ then $E(X)=\sum_{k=0}^n{n \choose k}P^k (1-P)^{n-k}=nP. $
\item If $X\sim $ Poisson distribution\index{Poisson distribution} ($\lambda$), then $E(X)= \sum_{k=0}^\infty k e^{-\lambda} \frac{\lambda^k}{k!}.$
\item If $X\sim$ Normal distribution\index{normal distribution} $(0,1)$ and $P(X=k)=0,$ then \begin{align*}
                                                                   E(X)&= \int_{-\infty}^\infty x dF(x)\\&= \int_{-\infty}^\infty x \frac{1}{\sqrt{2\pi}} e^{-\frac{1}{2}x^2}dx. 
                                                                  \end{align*}
\item If $X\sim $ Exponential distribution\index{exponential distribution} (1), and $F(x)=1-e^{-x},~~x \geq 0\Rightarrow dF(x)=e^{-x}dx,$ then $E(X)=\int_{-\infty}^\infty x dF(x)=\int^\infty_0 x e^{-x} dx.$
\item If $X$ is a continuous random variable with density $f$, that is $F(x)=\int_{-\infty}^xf(t)dt, $ so $dF(x)= f(x)dx$, then  $E(X)=\int_{-\infty}^\infty x dF(x)=\int_{-\infty}^\infty x f(x)dx.$  
\end{enumerate}

\chapter{Inequalities involving probabilites}

The study of inequalities in probability theory is very essential as they sometimes form the crux of the problem at hand that is to be solved. Ineuqualities inherently are quite beautiful and they never lose that beauty even if we put them in a probabilistic setting. This unit shall examine some basic and famous inequalities involving probabilities.

We begin with the following well known inequality and some of its corollaries.

\begin{thm}[Markov\index{Markov inequality}]
If $X$ is a non-negative valued
 random variable on a probability space $(\Omega,\mathcal{F},P)$, then $$P(X>a) \leq \frac{E(X)}{a}$$ for any positive number $a.$
\end{thm}

We have the following corollary of the above.

\begin{cor}
 If $X$ is a 
 random variable on a probability space $(\Omega,\mathcal{F},P)$, then $$P(|X|>a) \leq \frac{E(|X|)}{a}$$ for any  number $a>0.$
\end{cor}

We can generalize Markov's inequality in the following way.

\begin{thm}
 If $X$ is a non-negative valued
 random variable on a probability space $(\Omega,\mathcal{F},P)$ and $g$ is a non-negative measurable function on $(\mathbb{R}, \mathcal {B})$, then $$P(g(X)>a) \leq \frac{E(g(X))}{a}$$ for any positive number $a.$
\end{thm}

We now come to a very easy application of Markov's inequality given below.

\begin{thm}[Chebyshev\index{Chebyshev's inequality}]
If $X$ is a non-negative random variable on a probability space $(\Omega,\mathcal{F},P)$ then we have $$P(\mid X-E(X)\mid>\epsilon)\leq \frac{Var(X)}{\epsilon^2}$$ for any positive $\epsilon$.
\end{thm}
Before proceeding further we define the following.

\begin{definition}[Convex function\index{convex function}]
A function $f$ is said to be convex if for $0<t<1$ we have $f(tx+(1-t)y)\leq tf(x)+(1-t)f(y)$, where $x$ and $y$ are real numbers.
\end{definition}

We have now the following very famous result.

\begin{thm}[Jensen\index{Jensen's inequality}]
If $X$ is a non-negative valued random variable on a probability space $(\Omega,\mathcal{F},P)$, then $$\psi(E(X))\leq E(\psi(X))$$ where $\psi$ is a convex function.
\end{thm}

After the above inequalities we now turn our discussion to inequalities which involve probabilities of sums.

\begin{thm}[Hoeffding\index{Hoeffding's inequality}]
If $X_1, X_2, \ldots, X_n$ are independent random variables and $a_i\leq X_i\leq b_i$ then for any $\epsilon >0$ we have $$P(\sum_{i=1}^{n}(X_i-E(X_i))>n\epsilon)\leq e^{\frac{-n^2\epsilon^2}{2\sum_{i=n}^{n}(b_i-a_i)^2}}.$$
\end{thm}

From the above we get the following corollary very easily.

\begin{cor}
Let $X_1, X_2, \ldots, X_n$ be i.i.d random variables with $a_i\leq X_i\leq b_i$ and probability one, then for any $\epsilon>0$ we have $$P(\mid \overline{X}-E(X_1)\mid n\epsilon)\leq e^{\frac{-n\epsilon^2}{2(b-a)^2}}$$ where $\overline{X}=\frac{1}{n}\sum_{i=1}^{n}X_i$. 
\end{cor}

\begin{thm}[Cauchy-Schwarz\index{Cauchy-Schwarz inequality}]
If $X$ is a random variable then we have $E(X^2)\leq (E(X))^2.$
\end{thm}

\begin{thm}[Normal Tail inequality\index{normal tail inequality}]
If $X$ follows $N(0,1)$ and for $\epsilon>0$ where $X$ is a random variable then we have $$P(X>\epsilon)\leq \frac{1}{\epsilon}\frac{1}{\sqrt{2\pi}}e^{\frac{-1}{2}\epsilon^2}.$$
\end{thm}

\begin{thm}
If $X$ follows $N(0,1)$ and for $\epsilon>0$ where $X$ is a random variable then we have $$P(\mid X\mid >\epsilon)\leq \frac{2}{\epsilon \sqrt{2\pi}} e^{\frac{-1}{2}\epsilon^2}.$$
\end{thm}

We now have the following corollary from the above two theorems.

\begin{cor}
If $X_1, X_2, \ldots, X_n$ are i.i.d random variables and $X_i$ follows $N(0,1)$ then $\overline{X}$ also follows $N(0,1)$ and hence $\sqrt{n}\overline{X}$ follows $N(0,1)$ and we have $$P(\mid X\mid >\epsilon)\leq \frac{2}{\epsilon \sqrt{2n\pi}} e^{\frac{-1}{2}n\epsilon^2}$$ where $\overline{X}=\frac{1}{n}\sum_{i=1}^{n}X_i$. 
\end{cor}

We have the following simple result from analysis.

\begin{prop}
If $g(u)=-\theta u+\log (1-\theta+\theta e^4)$ where $0<\theta<1$ then $g(0)=0, g\prime(0)=0$ and $g\prime \prime (u)\leq \frac{1}{4}$.
\end{prop}

Using the above we get the following.

\begin{lem}
If $X$ is a random variable and $E(X)=0$ and $P(a\leq X\leq b)=1$ then for any $s>0$ we have $$E(e^{sX})\leq e^{\frac{s^2(b-a)^2}{8}}.$$
\end{lem}

The above lemma gives us another variant of the Hoeffding's inequality as stated below.

\begin{thm}[Hoeffding\index{Hoeffding's inequality}]
If $X_1, X_2, \ldots, X_n$ are independent random variables and $a_i\leq X_i\leq b_i$ then for any $\epsilon >0$ we have $$P(\sum_{i=1}^{n}(X_i-E(X_i))>n\epsilon)\leq e^{\frac{-2n^2\epsilon^2}{\sum_{i=n}^{n}(b_i-a_i)^2}}.$$
\end{thm}

From the above we get the following corollary very easily.

\begin{cor}
Let $X_1, X_2, \ldots, X_n$ be i.i.d random variables with $a_i\leq X_i\leq b_i$ and probability one, then for any $\epsilon>0$ we have $$P(\mid \overline{X}-E(X_1)\mid n\epsilon)\leq 2e^{\frac{-2n\epsilon^2}{(b-a)^2}}$$ where $\overline{X}=\frac{1}{n}\sum_{i=1}^{n}X_i$. 
\end{cor}

The following is a consequence of the above inequality.

\begin{thm}
Let $\{X_n\}_{n=1, 2, \cdots}$ be a sequence of independent Bernoulli $(p)$ random variables then $$\frac{1}{n}\sum_{i=1}^{n}X_i\rightarrow p$$ almost surely.
\end{thm}

We now come to an inequality involving medians.

\begin{thm}[Levy\index{Levy's inequality}]
If $X_1, X_2, \ldots, X_n$ are i.i.d random variables then $$P\{\max_{1\leq j\leq n}(S_j-m(S_j-S_n))\leq \epsilon \}\leq 2P(S_n\geq \epsilon)$$ where $\epsilon>0$ and $S_j=\sum_{i=1}^{j}X_i$ and $m(Y)$ denotes the median of $Y$.
\end{thm}

In this final part of the report we now mention two very famous inequalities related to expectation.

\begin{definition}
We define the norm of the expectation of a random variable as $\parallel X\parallel_p=(E(\mid X\mid^p))^{\frac{1}{p}}$ where $p>0$.
\end{definition}

\begin{thm}[Holder Inequality]\index{Holder's inequality}
For any $1<p<\infty$ and $1<q<\infty$ such that $\frac{1}{p}+\frac{1}{q}=1$ we have $$ E(\mid XY\mid)\leq \parallel X\parallel_p \parallel Y\parallel_q$$ where $X$ and $Y$ are two random variables.
\end{thm}

\begin{proof}
Let us define the following, $$ U=\frac{X}{\parallel X\parallel_p},~~V=\frac{Y}{\parallel Y\parallel_q}.$$
Now we have $$E(|UV|)\leq \frac{E(|U|^p)}{p}+\frac{E(|V|^q)}{q}=\frac{1}{p}+\frac{1}{q}=1.$$
Also $$E(|U|^p)=E\left(\frac{|X|^p}{E(|X|^P)}\right)=1.$$
Similarly we have $E(|V|^p)=1.$

Again from the above definition we have \begin{align*}
\parallel X\parallel_p&=(E(\mid X\mid^p))^{\frac{1}{p}}\\
\Rightarrow E\left(\left|\frac{XY}{\parallel X\parallel_p \parallel Y\parallel_q}\right|\right)&\leq 1\\
\Rightarrow~~~~~~~~~~~~~~~~ E(|XY|)&\leq \parallel X\parallel_p \parallel Y\parallel_q.
\end{align*}
\end{proof}

\begin{thm}[Minkowski Inequality]\index{Minkowski's inequality}
We have for any two random variables $X$ and $Y$ and $1\leq p< \infty$ $$ \parallel X+Y\parallel_p \leq \parallel X\parallel_p+\parallel Y\parallel_p.$$
\end{thm}

\begin{proof}
For $p=1$, the proof is trivial, as $E(|X+Y|)\leq E(|X|)+E(|Y|),~~i.e.,~~\\\parallel X+Y\parallel_1 \leq \parallel X\parallel_1+\parallel Y\parallel_1.$\\

For $1<p,~f(x)=x^p$ is a convex function on $(0,\infty)$ as $f'(x)\uparrow$. 

Now we define like earlier $$ U=\frac{X}{\parallel X\parallel_p},~~V=\frac{Y}{\parallel Y\parallel_q}.$$

Let $\lambda=\frac{\parallel X\parallel_p}{\parallel X\parallel_p+\parallel Y\parallel_p}$, then we have $1-\lambda=\frac{\parallel Y\parallel_p}{\parallel X\parallel_p+\parallel Y\parallel_p}.$

Now \begin{align*}
     E(|X+Y|^p)&\leq E(|X|+|Y|)^p,~~~\text{$p>1$}\\
&\leq E(\parallel X\parallel_p U+\parallel Y\parallel_p V)^p\\
&=E\{(\parallel X\parallel_p+\parallel Y\parallel_p)^p (\lambda U+(1-\lambda)V)^p\}\\
&=(\parallel X\parallel_p+\parallel Y\parallel_p)^p E(\lambda U+(1-\lambda)V)^p)\\
&\leq(\parallel X\parallel_p+\parallel Y\parallel_p)^p \{\lambda E(U^p)+(1-\lambda)E(V^p)\}.
    \end{align*}
Note that $E(U^p)=1$ and $E(V^p)=1$. So from above we have
\begin{align*}
E(|X+Y|^p)&\leq (\parallel X\parallel_p+\parallel Y\parallel_p)^p\\
\Rightarrow\parallel X+Y\parallel_p &\leq \parallel X\parallel_p+\parallel Y\parallel_p.
\end{align*}
\end{proof}



\fancyhead[LO]{}
\fancyhead[LE]{}
\fancyhead[RE]{}
\fancyhead[RO]{}
\fancyhead[LE]{\sffamily\leftmark}
\fancyhead[LO]{\sffamily\leftmark}


\clearpage
\fancyhead[LO]{}
\fancyhead[LE]{}
\fancyhead[RE]{}
\fancyhead[RO]{}
\fancyhead[LO]{}

\printindex

\end{document}